\documentclass[12pt]{amsart}
\usepackage{amssymb,amsmath}
\numberwithin{equation}{section}
\def\T{\text}
\def\simgeq{\underset\sim>}
\def\simleq{\underset\sim<}
\def\1#1{\overline{#1}}
\def\2#1{\widetilde{#1}}
\def\3#1{\widehat{#1}}
\def\4#1{\mathbb{#1}}

\def\5#1{\frak{#1}}
\def\6#1{{\mathcal{#1}}}
\def\C{{\4C}}
\def\R{{\4R}}

\def\Z{{\4Z}}

\begin{document}
\abstract
For the $\bar\partial$-Neumann problem on a regular coordinate domain $\Omega\subset \C^{n+1}$, we prove $\epsilon$-subelliptic estimates for  an index $\epsilon$ which is in some cases better  than $\epsilon=\frac1{2m}$ ($m$ being the {\it multiplicity}) as it was previously proved by Catlin and Cho in \cite{CC08}. This also supplies a much simplified proof of the existing literature. Our approach is founded on the method by Catlin in \cite{C87} which consists in constructing a family of weights $\{\phi^\delta\}$ whose Levi form is bigger than $\delta^{-2\epsilon}$ on the $\delta$-strip around $\partial\Omega$.
\newline
MSC: 32F10, 32F20, 32N15, 32T25 
\endabstract
\title[Precise subelliptic estimates for a class of special domains]{Precise subelliptic estimates for a class of special domains}
\author[T.V.~Khanh and G.~Zampieri]{Tran Vu Khanh and  Giuseppe Zampieri}
\address{Dipartimento di Matematica, Universit\`a di Padova, via 
Belzoni 7, 35131 Padova, Italy}
\email{khanh@math.unipd.it,  
zampieri@math.unipd.it}
\maketitle
\def\Giialpha{\mathcal G^{i,i\alpha}}
\def\cn{{\C^n}}
\def\cnn{{\C^{n'}}}
\def\ocn{\2{\C^n}}
\def\ocnn{\2{\C^{n'}}}
\def\const{{\rm const}}
\def\rk{{\rm rank\,}}
\def\id{{\sf id}}
\def\aut{{\sf aut}}
\def\Aut{{\sf Aut}}
\def\CR{{\rm CR}}
\def\GL{{\sf GL}}
\def\Re{{\sf Re}\,}
\def\Im{{\sf Im}\,}
\def\codim{{\rm codim}}
\def\crd{\dim_{{\rm CR}}}
\def\crc{{\rm codim_{CR}}}
\def\phi{\varphi}
\def\eps{\varepsilon}
\def\d{\partial}
\def\a{\alpha}
\def\b{\beta}
\def\g{\gamma}
\def\G{\Gamma}
\def\D{\Delta}
\def\Om{\Omega}
\def\k{\kappa}
\def\l{\lambda}
\def\L{\Lambda}
\def\z{{\bar z}}
\def\w{{\bar w}}
\def\Z{{\1Z}}
\def\t{{\tau}}
\def\th{\theta}
\emergencystretch15pt
\frenchspacing
\newtheorem{Thm}{Theorem}[section]
\newtheorem{Cor}[Thm]{Corollary}
\newtheorem{Pro}[Thm]{Proposition}
\newtheorem{Lem}[Thm]{Lemma}
\theoremstyle{definition}\newtheorem{Def}[Thm]{Definition}
\theoremstyle{remark}
\newtheorem{Rem}[Thm]{Remark}
\newtheorem{Exa}[Thm]{Example}
\newtheorem{Exs}[Thm]{Examples}
\def\Label#1{\label{#1}}
\def\bl{\begin{Lem}}
\def\el{\end{Lem}}
\def\bp{\begin{Pro}}
\def\ep{\end{Pro}}
\def\bt{\begin{Thm}}
\def\et{\end{Thm}}
\def\bc{\begin{Cor}}
\def\ec{\end{Cor}}
\def\bd{\begin{Def}}
\def\ed{\end{Def}}
\def\br{\begin{Rem}}
\def\er{\end{Rem}}
\def\be{\begin{Exa}}
\def\ee{\end{Exa}}
\def\bpf{\begin{proof}}
\def\epf{\end{proof}}
\def\ben{\begin{enumerate}}
\def\een{\end{enumerate}}
\def\dotgamma{\Gamma}
\def\dothatgamma{ {\hat\Gamma}}

\def\simto{\overset\sim\to\to}
\def\1alpha{[\frac1\alpha]}
\def\T{\text}
\def\R{{\Bbb R}}
\def\I{{\Bbb I}}
\def\C{{\Bbb C}}
\def\Z{{\Bbb Z}}
\def\Fialpha{{\mathcal F^{i,\alpha}}}
\def\Fiialpha{{\mathcal F^{i,i\alpha}}}
\def\Figamma{{\mathcal F^{i,\gamma}}}
\def\Real{\Re}
%
%
%
\section{  Introduction}
{\it Regular coordinate} domains  have $\epsilon$-subelliptic estimates; we discuss here about the optimal $\epsilon$.  These domains are defined by
\begin{equation}
\Label{1.1}
2\Re z_{n+1}+\underset{j=1}{\overset N\sum}|f_j(z)|^2<0,\quad \T{for }(z,z_{n+1})\in\C^n\times\C,
\end{equation}
where the $f_j$'s are holomorphic functions in a neighborhood of $0$
which satisfy $f_j=f(z_1,...,z_j)$ and $f_j(0,...,0,z_j)\neq0$. We denote by $m_j$ the smallest index such that $\partial_{z_j}^{m_j}f_j\neq0$ and define $m:=\underset{j=1,...,n}\Pi m_j$.
The  D'Angelo's type $D$ is the order of contact  of $\partial\Omega$ with any complex $1$-dimensional complex variety.
For these domains it is readily seen that $\sum_j|f_j|^2\simgeq |z|^{2m}$ which implies $D\leq 2m$.
 Moreover,  according to Catlin \cite{C83}, we must have $\epsilon\leq \frac1 D$ and, conversely, it is a conjecture by D'Angelo's \cite{D92} that $\epsilon\geq \frac1{2m}$. 
(Indeed, the conjecture is formulated for more general special domains and in this case the integer $m$ is the {\it multiplicity}.)
We present a simplified proof of a recent result by Catlin and Cho \cite{CC08} which gives positive answer to the conjecture for {  coordinate domains}. More important,  we find an intermediate number $\frac1{2m}\leq \frac\gamma 2\leq \frac1 D$ obtained by combining vanishing orders of the $f_j$'s in different directions and prove subelliptic estimates for $\epsilon$ coinciding with this new number $\frac\gamma2$. For instance, consider the domain defined by \eqref{1.1} for the choice $f_1=z_1^{m_1}$ and then
\begin{equation}
\Label{1.1bis}
f_j(z)=z_j^{m_j}+z_{j-1}^{l_j},\,\,l_j\leq m_{j-1},\,\,j=2,...,n.
\end{equation}
According to  Theorem~\ref{t2.2}, we set $\gamma_1=\frac1{m_1}$, define inductively $\gamma_j=\underset{i\leq j-1}\min\frac{l_i}{m_j}\gamma_i$ and write $\gamma$ for $\gamma_n$. 
(In particular, when $l_j= m_{j-1}$ and $m_j\geq m_{j-1}$ for any $j$, then $\gamma_j=\underset{i\leq j}\min\frac1{m_i}$.)
Then we have $\epsilon$-subelliptic estimates for $\epsilon=\frac\gamma2$. 

Note that we always have $\frac\gamma2\geq \frac1{2m}$ with equality holding only if $l_j=1$ for any $j$. 
On the other hand we have
\begin{equation}
\Label{1.2}
\frac\gamma2\geq\frac1 D,
\end{equation}
which shows that our index of subellipticity is optimal in this example.
To prove \eqref{1.2}, we have just to notice that  the curve
$$
\C\to\C^{n+1},\quad\tau\mapsto(\tau^{\frac{\gamma_1}\gamma},...,\tau^{\frac{\gamma_{n-1}}\gamma},\tau,0),
$$
has order of contact equal to $\frac2\gamma$.
On the other hand, the presence of $\epsilon$-subelliptic estimates for $\epsilon=\frac\gamma2$ assured by our theorem, yields equality in \eqref{1.2}.

\section{Precise subellipticity index for a class of special domains}
Let $z=(z_1,...,z_n)$ be coordinates in $\C^n$ and $(z,z_{n+1})$ coordinates in $\C^n\times\C$. We deal with { regular coordinate} domains $\Omega\subset\C^{n+1}$, that is, domains defined by \eqref{1.1} for $f_j$ holomorphic which satisfy $f_j=f_j(z_1,...,z_j)$ and $\partial^{m_j}_{z_j}f_j\neq0$ for some $m_j$.  For these domains we consider the $\bar\partial$-Neumann problem and, in particular, the {\it $\epsilon$-subelliptic estimates}. These are of the type
\begin{equation}
\Label{2.3}
|||u|||^2_\epsilon\simleq ||\bar\partial u||^2_0+||\bar\partial^*u||^2_0+||u||^2_0,
\end{equation}
for any $C^\infty_c(\bar\Omega)$-form $u$ of degree $k\geq 1$ in the domain $D_{\bar\partial^*}$ of $\bar\partial^*$. Here $|||\cdot|||_\epsilon$ is the {\it tangential} Sobolev norm of index $\epsilon$. This is said an estimate with $\epsilon$-fractional gain of derivative. It is classical (see \cite{FK72}) that it implies the local hypoellipticity of the $\bar\partial$-Neumann problem. The canonical solution $u$ of $\bar\partial u=f$, that is, the solution orthogonal to $\ker\bar\partial$, is $C^\infty$ up to $\partial\Omega$ precisely at those points of $\partial\Omega$ where $f$ is $C^\infty$. In particular, the Bergman projection preserves $C^\infty$ smoothness. 
The following theorem has been  recently obtained by \cite{CC08}; we give here a much simplified proof.
\bt
\Label{t2.1}
Let $\Omega\subset\C^{n+1}$ be a regular coordinate domain defined by \eqref{1.1} with the $f_j$'s satisfying $f_j=f_j(z_1,...,z_j)$. Write 
\begin{equation}
\Label{2.4}
f_j=z_j^{m_j}+O(z_1,...,z_{j-1},z_j^{m_j+1});
\end{equation}
then $\epsilon$-subelliptic estimates hold for
$\epsilon=\frac1{2m_1\cdot...\cdot m_n}.$
\et
\bpf
According to Catlin \cite{C87} it suffices to find a family of bounded weights $\{\phi^\delta\}$ for $\delta\searrow0$ whose Levi form satisfies $\sum_{ij}\phi^\delta_{ij}u_i\bar u_j\geq \delta^{-\frac1{m_1\cdot...\cdot m_n}}|u|^2$ over the strip of $\Omega$ about the boundary $S_\delta=\{z\in\Omega:\,-r(z)<\delta\}$. 
Once the functions $\phi^\delta$ have been found, we have to deform them to new functions $\tilde \phi^\delta$, bounded and plurisubharmonic not only in $S_\delta$ but in the whole $\Omega$, satisfying the same Levi conditions as the $\phi^\delta$  on $S_{\frac\delta2}$; for this we refer to Lemma~\ref{l2.1} after the end of the proof of Theorem~\ref{t2.1}.
To define the functions $\phi^\delta$, we put $\gamma_0=1$, $\gamma_j=\frac1{m_1\cdot...\cdot m_j}$, choose $\alpha\geq 1$ and put $\alpha_j=\alpha m_{j+1}...m_n$. Note that $(\alpha_{j-1}-1)-\alpha_j(m_j-1)=\alpha_j-1.$ We choose a cut off function $\chi$ which satisfies
\begin{equation*}
\begin{cases}
\T{supp}\,\chi\subset(0,2),
\\
\chi\equiv1\T{ on }(0,1),
\end{cases}
\end{equation*}
and also take a small constant $c$ to be specified later.
We also rewrite the inequality \eqref{1.1} which defines $\Omega$ as   $r<0$ and
define
\begin{equation}
\Label{10}
\begin{split}
\phi^\delta=-\log(\frac{-r+\delta}\delta)+&\underset{j=1}{\overset n\sum}\underset{h=1}{\overset{m_j-1}\sum}\frac1{|\log\,*|}\log\left(|\partial_{z_j}^hf_j|^2+\frac{\delta^{(m_j-h)\gamma_j}}{|\log\delta|^{(m_j-h)\alpha_j}}\right)
\\
&+c\underset{j=1}{\overset n\sum}\chi(\frac{|z_j|^2}{\delta^{\gamma_j}})\log\left(\frac{|z_j|^2+\delta^{\gamma_j}}{\delta^{\gamma_j}}\right),
\end{split}
\end{equation}
where $*=\frac{\delta^{\gamma_j(m_j-h)}}{|\log\delta|^{(m_j-h)\alpha_j}}$. Notice here that $\log *\simeq \log \delta$. The weights $\phi^\delta$ that we have defined are bounded in the strip $S_\delta$. 
We use the notations
\begin{equation*}
A_j=\delta^{-1}|\partial f_j\cdot u|^2\qquad B_j^h=\delta^{-(m_j-h)\gamma_j}|\partial\partial_{z_j}f_j\cdot u|^2.
\end{equation*}
We also denote by $ cC_j$ the Levi form of the third term of \eqref{10} applied to $u$; note that $C_j\simgeq 
\delta^{-\gamma_j}|u_j|^2$ for $|z_j|\leq \delta^{\gamma_j}$, but, otherwise, $ C_j$ can take negative values; however, $|C_j|\simleq \delta^{-\gamma_j}|u_j|^2$. We also set
$$
 D_j=A_j+\underset{h\leq k_j}\sum B_j^h+c C_j.
$$
We first prove an auxiliary  statement: the assumption
\begin{equation}
\Label{20}
\underset{i\leq i_o}\sum D_i\geq \delta^{-\gamma_{i_o}}|\log\delta|^{\alpha_{i_o}-1}|u_{i_o}|^2+\underset{i\leq i_o-1}\sum\delta^{-\gamma_i}|u_i|^2,
\end{equation}
implies, for any $j\geq i_o+1$
\begin{equation}
\Label{21}
\underset{i\leq j}\sum D_i\geq \delta^{-\gamma_j}|\log\delta|^{\alpha_j-1}\underset{i\leq j}\sum|u_i|^2.
\end{equation}
In fact, 
by iteration, it suffices to prove the statement for $i_o=j-1$. For that,
 we first notice that $\delta^{-\gamma_i}\geq \delta^{-\gamma_{i_o}}|\log\delta|^k$ for any $k$ and  for any $i\leq i_o-1$. It follows
\begin{equation}
\Label{14}
\begin{split}
A_j+\underset{i\leq j-1}\sum D_i&\simgeq \delta^{-\gamma_{j-1}}|\log|^{\alpha_{j-1}-1}\left(|\partial_{z_j}f_j|^2|u_j|^2-\underset{i\leq j-1}\sum|u_i|^2\right)+\underset{i\leq j-1}\sum D_i
\\
&\simgeq \delta^{-\gamma_{j-1}}|\log|^{\alpha_{j-1}-1}|\partial_{z_j}f_j|^2|u_j|^2.
\end{split}
\end{equation}
Assume at this point $|\partial_{z_j}f_j|\geq \frac{\delta^{(m_j-1)\gamma_j}}{|\log\delta|^{(m_j-1)\alpha_j}}$; then \eqref{14} can be continued by
\begin{equation*}
\begin{split}
{}&
\geq \delta^{-\gamma_{j-1}+(m_j-1)\gamma_j}|\log\delta|^{(\alpha_{j-1}-1)-(m_j-1)\alpha_j}
\\
&\geq \delta^{-\gamma_j}|\log\delta|^{\alpha_j-1}|u_j|^2.
\end{split}
\end{equation*}
Notice that  this controls $c C_j$ when this gets negative; this happens in all cases which follow; hence we avoid to recall it at each step. If not, we can assume $|\partial_{z_j}f_j|\leq \frac {\delta^{(m_j-1)\gamma_j}}{|\log \delta|^{(m_j-1)\alpha_j}}$ and then use $B^1_j$. We have
\begin{equation}
\Label{15}
\begin{split} 
B_j^1+\underset{i\leq j-1}\sum D_i&\simgeq \delta^{-\gamma_j(m_j-1)}\frac{|\log\delta|^{(m_j-1)\alpha_j}}{|\log *|}
\left(|\partial_{z_j}^2f_j|^2|u_j|^2-\underset{i\leq j-1}\sum|u_i|^2\right)+\underset{i\leq j-1}\sum D_i
\\
&\simgeq \delta^{-\gamma_{j-1}+\gamma_j}|\log\delta|^{(m_j-1)\alpha_j-1}|\partial_{z_j}^2f_j|^2|u_j|^2.
\end{split}
\end{equation}
If $|\partial_{z_j}^2 f_j|\geq \frac{\delta^{(m_j-2)\gamma_j}}{|\log\delta|^{(m_j-2)\alpha_j}}$, then
\eqref{15} can be continued by 
$$
\geq \delta^{-\gamma_{j-1}+\gamma_j+\gamma_{j-1}-2\gamma_j}|\log\delta|^{\alpha_j-1}|u_j|^2.
$$
If $|\partial_{z_j}^2f_j|\leq \frac{\delta^{(m_j-2)\gamma_j}}{|\log\delta|^{(m_j-2)\alpha_j}}$, we pass to $B^2_j$. In this way we jump from $\partial^{h-1}_{z_j}f_j$ to $\partial^{h}_{z_j}f_j$ until we reach $B_j^{m_j-1}$. At this stage we have $|\partial_{z_j}^{m_j-1}f_j|\leq \frac {\delta^{\gamma_j}}{|\log\delta|^{\alpha_j-1}}$ which yields readily
$$
B_j^{m_j-1}\geq \delta^{-\gamma_j}|\log \delta|^{\alpha_j-1}|u_j|^2.
$$
This proves that \eqref{20} for $i_o=j-1$  implies \eqref{21}. By iteration the conclusion is true for any $i_o\leq j-1$. 
We now prove that, for any value of $z_j$
\begin{equation}
\Label{11}
\underset{i\leq j}\sum D_i\simgeq \delta^{-\gamma_j}|u_j|^2,
\end{equation}
whereas, when $|z_j|\geq \delta^{\gamma_j}$
\begin{equation}
\Label{12}
\underset{i\leq j}\sum D_i\simgeq
\begin{cases}
\T{\rm either} &\delta^{-\gamma_j}|\log\delta|^{\alpha_j-1}|u_j|^2
\\
\T{\rm or}&\delta^{-\gamma_{j-1}}|z_j|^{2(m_j-1)}|u_j|^2.
\end{cases}
\end{equation}

We proceed by induction over $j$. The first step $j=1$ is easy. In fact
$ C_1\simgeq \delta^{-\gamma_1}|u_1|^2$ for $|z_1|\leq \delta^{\gamma_1}$ and
$$
A_1\simgeq \max(\delta^{-\gamma_1},\delta^{-1}|z_1|^{2(m_1-1)})|u_1|^2\quad\T{ for $|z_1|\geq \delta^{\gamma_1}$}.
$$
In particular, when $cC_j$ takes negative values, these are controlled by $A_1$ for suitably small $c$.
This proves \eqref{11} and \eqref{12} for $j=1$. Suppose that \eqref{11} and \eqref{12} are true up to step $j-1$ and prove them for $j$.  If $|z_j|\leq \delta^{\gamma_j}$, then $ C_j\simgeq \delta^{-\gamma_j}|u_j|^2$. Otherwise, assume $|z_j|\geq \delta^{\gamma_j}$. 
We show now that, under this assumption, we must have \eqref{20} for some $i_o\leq j-1$ unless the second alternative in \eqref{12} holds. In fact,
let $|z_{j-1}|\simgeq  |z_j|^{m_j-1}$; then $|z_{j-1}|\simgeq \delta^{\gamma_j(m_j-1)}\geq \delta^{\gamma_{j-1}-\gamma_j}$ and therefore 
\begin{equation*}
\begin{split}
\delta^{-\gamma_{j-2}}|z_{j-1}|^{2(m_{j-1}-1)}
&\geq \delta^{-\gamma_{j-2}+\gamma_j(m_{j-1}-1)(m_j-1)}
\\
&=\delta^{-\gamma_{j-2}+\gamma_{j-2}-\gamma_{j-1}-\gamma_jm_{j-1}+\gamma_j}
\\
&\geq \delta^{-\gamma_{j-1}}|\log\delta|^k\quad\T{for any choice of $k$.}
\end{split}
\end{equation*}
On the other hand, for any choice of $i\leq j-2$, we have 
$$
\delta^{-\gamma_i}\geq \delta^{-\gamma_{j-1}}|\log\delta|^k.
$$
This implies \eqref{20} for $i_o=j-1$ which implies in turn \eqref{21}. Otherwise, we assume $|z_{j-1}|\simleq |z_j|^{m_j-1}$. Now, we point our attention to $z_{j-2}$. If $|z_{j-2}|\simgeq |z_j|^{m_j-1}$, then $|z_{j-2}|\geq \delta^{\gamma_j(m_j-1)}\geq \delta^{\gamma_{j-1}-\gamma_j}$ and hence
\begin{equation*}
\begin{split}
\delta^{-\gamma_{j-3}}|z_{j-2}|^{2(m_{j-2}-1)}&\geq \delta^{-\gamma_{j-3}+(m_{j-2}-1)(\gamma_{j-1}-\gamma_j)}
\\
&=\delta^{-\gamma_{j-3}+m_{j-2}\gamma_{j-1}}
\\
&\geq\delta^{-\gamma_{j-2}}|\log\delta|^k.
\end{split}
\end{equation*}
On the other hand, we have for any $i\leq j-3$, $\delta^{-\gamma_i}\geq \delta^{-\gamma_{j-2}}|\log\delta|^k$. This implies
\begin{equation}
\Label{22}
\underset{i\leq j-2}\sum D_i\geq \delta^{-\gamma_{j-2}}|\log\delta|^{\alpha_{j-2}-1}\underset{i\leq j-2}\sum|u_i|^2.
\end{equation}
The same argument which shows that \eqref{20} implies \eqref{21} also serves in proving that \eqref{22} implies \eqref{20}; in turn, \eqref{20} implies \eqref{21}. We can therefore assume
$$
|z_{j-2}|\simleq |z_j|^{m_j-1}.
$$
We repeat the argument that we developed for $i=j-1$ and $i=j-2$ for any other index $i\leq j-1$. 
This can be explained by the fact that the second of \eqref{12} implies the first at each step $i\leq j-1$ (with the cases $i=j-1$ and $i=j-2$ having already been proved). In fact, if $|z_i|\geq |z_j|^{m_j-1}$ and hence $|z_i|\geq \delta^{\gamma_{j-1}-\gamma_j}$, then 
\begin{equation*}
\begin{split}
\delta^{-\gamma_{i-1}}|z_i|^{2(m_i-1)}&\geq \delta^{-\gamma_{i-1}+(\gamma_{j-1}-\gamma_j)(m_i-1)}
\\
&\geq \delta^{-\gamma_{i-1}+m_i\gamma_{j-1}}
\\
&\geq \delta^{-\gamma_i}|\log\delta|^k.
\end{split}
\end{equation*}
This yields \eqref{20} and thus also \eqref{21} unless
\begin{equation}
\Label{23}
|z_i|\leq |z_j|^{m_j-1}\qquad\T{\rm for any $i\leq j-1$}.
\end{equation}
On the other hand, when \eqref{23} holds, then
\begin{equation*}
|\partial_{z_j}f_j|\simgeq |z_j|^{2(m_j-1)}-\frac12\underset{i\leq j-1}\sum |z_i|^2\simgeq |z_j|^{2(m_j-1)}.
\end{equation*}
It follows
\begin{equation}
\Label{24}
\begin{split}
A_j+\underset{i\leq j-1}\sum D_i&\geq \delta^{-1}|\partial f_j\cdot u|^2+\underset{i\leq j-1}\sum D_i
\\
&\simgeq \delta^{-\gamma_{j-1}}(|z_j|^{2(m_j-1)}|u_j|^2-\underset{i\leq j-1}\sum|u_i|^2)+\underset{i\leq j-1}\sum D_i
\\
&\simgeq \delta^{-\gamma_{j-1}}|z_j|^{2(m_j-1)}|u_j|^2.
\end{split}
\end{equation}
Note that this is in any case $\simgeq \delta^{-\gamma_j}|u_j|^2$ and also that it controls $c\tilde C_j$, for suitable $c$ when $\tilde C_j$ gets negative. So, in this case we have the second alternative in \eqref{12}.      
This concludes the proof of the theorem.

\epf
It remains to prove the technical lemma which shows how to modify the functions $\phi^\delta$ to $\tilde\phi^\delta$ so that they are plurisubharmonic on the whole of $\Omega$.
\bl
\Label{l2.1}
There are $\tilde\phi^\delta$, plurisubharmonic and bounded on $\Omega$ and such that
\begin{equation*}
\tilde \phi^\delta=
\begin{cases}
e^{2\frac r\delta}+e^{-2}\phi^\delta&\T{ on $S_{\frac\delta2}$},
\\
0&\T{ on $\Omega\setminus S_\delta$}.
\end{cases}
\end{equation*}
\el
\bpf
We take $\theta:\R\to\R^+,\,\,t\mapsto\theta(t)$, convex increasing, that is, satisfying $\dot\theta\geq0,\,\ddot\theta\geq0$ and such that
\begin{equation*}
\theta=
\begin{cases}0&\T{ for $t\leq 2e^{-2}$},
\\
t&\T{ for $t\geq e^{-1}$},
\end{cases}
\end{equation*}
set $\psi^\delta:=e^{2\frac r\delta}+e^{-2}\phi^\delta$ and define $\tilde \phi^\delta:=\theta\circ\psi^\delta$. Remember that $\phi^\delta$ take values in $[0,1]$ and notice that $S_{\frac\delta2}\subset\{z:\,\psi^\delta\geq e^{-1}\}$ and $\Omega\setminus S_\delta\subset\{z:\,\psi^\delta\leq 2e^{-2}\}$; thus $\tilde \phi^\delta=\psi^\delta$ on $S_{\frac\delta2}$ and $\tilde\phi^\delta=0$ on $\Omega\setminus S_\delta$.

\epf

We now prove subelliptic estimates for a class of domains with a better index $\epsilon\geq \frac 1{2m_1\cdot...\cdot m_n}$; in particular, for the domains defined by \eqref{1.1}, this index coincides with the optimal value $\epsilon=\frac1D$. To achieve our goal, we have to specify the vanishing order of $f_j$ in different directions $z_i$ for $i\leq j$.
\bt
\Label{t2.2}
Let $\Omega\subset\C^{n+1}$ be a regular coordinate domain defined by \eqref{1.1} with $f_j=f_j(z_1,...,z_j)$ and $\partial^{m_j}_{z_j}f_j\neq0$. We  denote by $l_j^i,\,\,i<j,$ the vanishing order of $f_j$ in $z_i$, that is, we write
\begin{equation}
\Label{2.10}
f_j=z_j^{m_j}+O_j(z_1^{l_j^1},...z_{j-1}^{l_j^{j-1}},z_j^{m_j})\quad\T{\rm for $l_j^i\leq m_i$}.
\end{equation}
Assume that each $O_j$ contains no power of $z_j$ of degree $\leq m_j-1$ and define $\gamma_1=\frac1{m_1}$ and, inductively,
$$
\gamma_j=\underset{i\leq j-1}\min\frac{l_j^i}{m_j}\gamma_i
$$
and also write $\gamma$ for $\gamma_n$. 
Then we have $\epsilon$-subelliptic estimates for any $\epsilon=\frac\gamma2$.
\et
\bpf
Let $k_j$ be the highest power $\leq m_j-1$ of $z_j$ in $O_j$
and let $r<0$ be the inequality which defines $\Omega$.
 We define 
\begin{equation}
\Label{2.11}
\phi^\delta=-\log(\frac{-r+\delta}\delta)+c\underset{j=1}{\overset n\sum}\chi(\frac{|z_j|^2}{\delta^{\gamma_j}})\log\left(\frac{|z_j|^2}{\delta^{\gamma_j}}+1\right)
\end{equation}
and prove that for any $u\in C^n$,  $\phi^\delta$ satisfy $\sum_{ij}\phi^\delta_{ij}u_i\bar u_j\geq\delta^{-\gamma}|u|^2$ over the strip $S_\delta$. We denote by $c\underset{j=1}{\overset n\sum}\phi_j$ the second term in the right hand side of \eqref{2.11} and define
\begin{gather*}
A_j=\delta^{-1}|\partial f_j\cdot u|^2,\qquad 
C_j=\partial\bar\partial \phi_j(u,\bar u) ,\qquad D_j=A_j+ cC_j.
\end{gather*}
We wish to prove that
\begin{equation}
\Label{2.12}
\underset{i\leq j}\sum D_i\simgeq
\underset{i\leq j}\sum
\delta^{-s_i\gamma_i}|z_i|^{2(s_i-1)}|u_i|^2\T{ for any $s_i\leq m_i$}.
\end{equation}
 It is easy to prove the first step, that is,
$$
D_1\simgeq \delta^{-s\gamma_1}|z_1|^{2(s-1)}|u_1|^2\quad\T{for any $s\leq m_1$}.
$$
Suppose we have already proved \eqref{2.12} for any $i\leq j-1$. We have to  prove that
\begin{equation}
\Label{2.8bis}
\underset{i\leq j}\sum D_i\geq \delta^{-s\gamma_j}|z_j|^{2(s-1)}|u_j|^2\quad\T{for any $s\leq m_j$}.
\end{equation}
We recall here that $l_j^i\leq m_i$ and notice that $\gamma_i\leq \frac1{m_i}$; in particular, $l_j^i\gamma_i\leq1$. We have
\begin{equation*}
A_j+\underset{i\leq j-1}\sum D_i\simgeq \delta^{-1}\left| \partial f_j\cdot u\right|^2
+\underset{i\leq j-1}\sum\delta^{-s_i\gamma_i}|z_i|^{2(s_i-1)}|u_i|^2.
\end{equation*}
We fix our choice of the $s_i$'s as $s_i=l_j^i$ which are smaller than $m_i$; thus $l_j^i\gamma_i\leq 1$. 
It follows
$$
A_j+\underset{i\leq j-1}\sum D_i\geq \underset{i\leq j-1}\sum\delta^{-l_j^i\gamma_i}\left[|\partial f_j\cdot u|^2+ |z_i|^{2(l_j^i-1)}|u_i|^2\right].
$$
On the other hand
$$
|\partial f_j\cdot u|^2\simgeq |z_j|^{2(m_j-1)}|u_j|^2-\underset{i\leq j-1}\sum |z_i|^{2(l_j^i-1)}|u_i|^2,
$$
and therefore, since $ l_j^i\gamma_i\geq m_j\gamma_j$ for any $i\leq j-1$, we conclude
\begin{equation}
\Label{30}
A_j+\underset{i\leq j-1}\sum D_i\simgeq\delta^{-m_j\gamma_j}|z_j|^{2(m_j-1)}|u_j|^2.
\end{equation}
This proves \eqref{2.12} for the choice $s=m_j$. 
We prove now \eqref{2.12} for $s=1$; for this we have to call into play $C_j$. We have
\begin{equation}
\Label{31}
\frac12C_j+\delta^{-m_j\gamma_j}|z_j|^{2(m_j-1)}\simgeq (\delta^{-\gamma_j}+\delta^{-m_j\gamma_j}|z_j|^{2(m_j-1)})|u_j|^2.
\end{equation}
In fact, if $|z_j|^2\leq \delta^{\gamma_j}$, then $C_j\geq \delta^{-\gamma_j}|u_j|^2$. If, instead, $|z_j|\geq \delta^{\gamma_j}$, and thus $C_j$ gets negative, we have on our side the fact that $\delta^{-m_j\gamma_j}|z_j|^{2(m_j-1)}|u_j|^2\geq \delta^{-\gamma_j}|u_j|^2$ and therefore it controls $cC_j$ for suitably small $c$. From \eqref{30} we conclude that $A_j+cC_j+\underset{i\leq j-1}\sum D_i$ is bigger than the right side of \eqref{31}, which yields \eqref{2.12} for $s=1$ and  $s=m_j$. The estimate \eqref{2.12} for general $s$ with $1\leq s\leq m_j$ is just a combination of the two opposite cases $s=1$ and $s=m_j$.

\epf

\br
The theorem applies in particular to the class of examples described by \eqref{1.1bis}.
\er
We have a final statement which collects  in a unified frame the conclusions of Theorems \ref{t2.1} and \ref{t2.2}.
\bt
\Label{t2.3}
 Let $\Omega\subset\C^{n+1}$ be a regular coordinate domain defined by \eqref{1.1} for $f_j=f_j(z_1,...,z_j)$ which satisfy \eqref{2.10}. 
Define $\gamma_1=\frac1{m_1}$ and, inductively
\begin{equation*}
 \gamma_j=
\begin{cases}
\underset{i\leq j-1}\min\frac1{m_j}\gamma_i&\T{\rm if $O_j$ contains some power of $z_j$ in degree $\leq m_j-1$},
\\
\underset{i\leq j-1}\min\frac{l^i_j}{m_j}\gamma_i&\T{\rm otherwise}
\end{cases}
\end{equation*}
and
also write $\gamma$ for $\gamma_n$.
Then we have $\epsilon$-subelliptic estiates for $\epsilon=\frac\gamma2$.
\et
\bpf
We use in this situation the family of weights
\begin{equation*}
\begin{split}
\phi^\delta=-\log(\frac{-r+\delta}\delta)&+
+\underset{\{j:k_j\neq0\}}{\sum}\underset{h\leq k_j}\sum\frac1{|\log*|}\log\left(|\partial_{z_j}^hf_j|^2+\frac{\delta^{(m_j-h)\gamma_j}}{|\log\delta|^{(m_j-h)\alpha_j}}\right)
\\
&+c\underset{j=1}{\overset n\sum}\chi\left(\frac{|z_j|^2}{\delta^{\gamma_j}}\right)
\log\left(\frac{|z_j|^2}{\delta^{\gamma_j}}+1\right).
\end{split}
\end{equation*}
The proof of he thorem is a combination of those of Theorem \ref{t2.1} and \ref{t2.2}.

\epf
\be
Let us consider in $\C^4$ the domain defined by
$$
2\Re  z_4+|z_1^6|^2+|z_2^4-z_1z_2|^2+|z_3^4-z_2^3+z_1|^2<0.
$$
Here $\gamma_1=\frac16$, $\gamma_2=\frac1{6\cdot 4}$ and $\gamma_3=\frac 3{6\cdot 4\cdot 4}$;  we have $\epsilon$-subelliptic estimates for $\epsilon=\frac{\gamma_3}2$.
\ee

\end{document}